\documentclass{amsart} 
\usepackage{amssymb}
\usepackage{amsmath}
\usepackage{amsfonts}

\begin{document}
\newtheorem{theo}{Theorem}[section]
\newtheorem{prop}[theo]{Proposition}
\newtheorem{lemma}[theo]{Lemma}
\newtheorem{coro}[theo]{Corollary}
\theoremstyle{definition}
\newtheorem{exam}[theo]{Example}
\newtheorem{defi}[theo]{Definition}
\newtheorem{rem}[theo]{Remark}


\newcommand{\Bb}{{\bf B}}
\newcommand{\Cb}{{\bf C}}
\newcommand{\Nb}{{\bf N}}
\newcommand{\Qb}{{\bf Q}}
\newcommand{\Rb}{{\bf R}}
\newcommand{\Zb}{{\bf Z}}
\newcommand{\Ac}{{\mathcal A}}
\newcommand{\Bc}{{\mathcal B}}
\newcommand{\Cc}{{\mathcal C}}
\newcommand{\Dc}{{\mathcal D}}
\newcommand{\Fc}{{\mathcal F}}
\newcommand{\Ic}{{\mathcal I}}
\newcommand{\Jc}{{\mathcal J}}
\newcommand{\Kc}{{\mathcal K}}
\newcommand{\Lc}{{\mathcal L}}
\newcommand{\Mx}{{\mathcal M}}
\newcommand{\Nc}{{\mathcal N}}
\newcommand{\Oc}{{\mathcal O}}
\newcommand{\Pc}{{\mathcal P}}
\newcommand{\Qc}{{\mathcal Q}}
\newcommand{\Sc}{{\mathcal S}}
\newcommand{\Tc}{{\mathcal T}}
\newcommand{\Uc}{{\mathcal U}}
\newcommand{\Vc}{{\mathcal V}}
\newcommand{\btu}{\bigtriangleup}

\author{Charles Akemann and Nik Weaver}

\title [Classically normal pure states]
       {Classically normal pure states}

\address {Department of Mathematics\\
          University of California\\
          Santa Barbara, CA 93106}
\address {Department of Mathematics\\
          Washington University in Saint Louis\\
          Saint Louis, MO 63130}

\email {akemann@math.ucsb.edu, nweaver@math.wustl.edu}

\date{May 6, 2007}

\begin{abstract}
 A pure state $f$ of a von 
Neumann algebra $\Mx$ is called {\it classically normal} if $f$ is normal on
any von Neumann subalgebra of $\Mx$ on which
$f$ is multiplicative.
Assuming the continuum hypothesis, a separably represented von
Neumann  algebra $M$ has classically normal, singular pure states iff there is a
central projection $p \in M$ such that $pMp$ is a factor of type $I_\infty, II$,
or
$III$. 
\end{abstract}

\maketitle

DEFINITION:   A pure state $f$ of a von 
Neumann algebra $\Mx$ is called {\it classically normal} if $f$ is normal on any
von Neumann subalgebra of $\Mx$ ("subalgebra" implies the same unit) on which $f$
is multiplicative.

By Lemma 0.2 below, a pure state $f$ on a von Neumann algebra $\Mx$
is  classically normal if, for every von
Neumann subalgebra $\Cc$ of $\Mx$, either $f$ is not multiplicative on $\Cc$, or
else there is a minimal projection $q$ in $\Cc$ such that $f(q)=1$ and $q$ is
central in $\Cc$.   Using the continuum hypothesis and a transfinite
construction, in Theorem 0.7 we show the existence of classically normal,
singular pure states on all infinite dimensional factors acting on a separable
Hilbert space. Corollary 0.8 contains the easy "only if" part of the main
result. 

 Here is some history. Let $H$ be a separable
infinite-dimensional Hilbert space and let
$\Bc(H)$ denote the algebra of all
bounded linear operators on $H$. Kadison and Singer \cite{ks} suggested that
every pure state on
$\Bc(H)$ would restrict to a pure state on some maximal abelian self-adjoint
subalgebra (aka MASA). Anderson \cite{j} formulated the stronger conjecture
that every pure state on $\Bc(H)$ is of the form
$f(a) = \lim_\Uc \langle ae_n, e_n\rangle$ for some orthonormal basis
$(e_n)$ and some ultrafilter $\Uc$ over the natural numbers $\Nb$.  Using the
continuum hypothesis, we showed in \cite{cn2} that these conjectures are false
by showing that there is a pure state $f$ on $\Bc(H)$ that is not
multiplicative on any MASA.  The argument in the key lemma of that paper used
powerful results about the Calkin algebra, so finding the "right" definitions
and proofs for general von Neumann algebras took some time.

Our construction of a classically normal pure state will be by
transfinite induction, just as in \cite{cn2}.  The
difference will be in the proofs of the Lemmas that allow the
transfinite construction to go through. We start with some easy facts.  

\begin{lemma}\label{LM}
Let  $f$ denote a state on a C*-algebra $\Bc$ in which the linear
combinations of the projections are dense.
$f$ is multiplicative on $\Bc$ iff $f(p) \in \{0,1\}$ for all projections
$p\in \Mx$. 
\end{lemma}
\begin{proof}
Suppose that $a,b \in B$ and
$f(ab)\ne f(ba)$  WLOG we can assume that 
$a=\sum s_jp_j, b = \sum t_iq_i$, finite
linear combinations of projections since the map
$(a,b) \to (ab-ba)\to f(ab-ba)$ is
continuous, so we only need to show that it
annihilates a dense set in $\Ac\times A$.
Then
$$f(ba) = \sum s_j t_jf(q_ip_j).$$
 $$f(ab) = \sum s_j t_jf(p_jq_i).$$

Thus it suffices to prove that for any projections $p,q \in B, f(pq)=f(qp)$. 
If either $f(p)=0$ or $f(q)=0$, then $0 \le
|f(pq)|=|f(qp)|\le f(p)^{1/2}f(q)^{1/2}=0$  by the Cauchy-Schwarz inequality. 
If
$f(q)=f(p)=1$, then $f(1-p)=f(1-q)=0$, so
$$f(pq)=f(q)-f((1-p)q)=1-0=1=f(p)=f(p)-f((1-q)p) = f(qp).$$
\end{proof}

\begin{lemma}\label{LM}Let $f$ denote a pure state on a von neumann algebra
$\Nc$.

 a.   $f$
is normal iff there is a minimal projection $p \in \Nc$ such that $f(p)=1$.  $f$
is both normal and multiplicative iff the projection $p$ (of the previous
sentence) is central.

b.  Let $K_f=\{x
\in \Mx$ such that $f(xx^*+x^*x)=0$\}.  Then $f$ is singular iff any
increasing, positive approximate unit of $K_f$ converges to 1 in $\Mx$ for the
weak* topology.
\end{lemma}
\begin{proof}
Let $f$ denote a normal pure state and $p$ is its support projection
(\cite{tak}, p. 140) in
$\Nc$.  Since the support projection $q$ of $f$ in $N^{**}$ is minimal in
$N^{**}$ (see sect. 3.13 of
\cite{gert}), $p$ must be minimal in $\Nc$ and $p=q$. 

Now suppose that $f$ is a pure state of $\Nc$ such that there is a minimal
projection $p$ in $\Nc$ such that $f(p)=1$. Then $f(a)=f(pap)$ by Cauchy-Shwarz
inequality.  Thus $f$ is normal because $p\Nc p$ is one dimensional (hence
$f|_{p\Nc p}$ is normal) and $a \to pap$ is weak* continuous. 

If $p$ is central, then $f(ab)=f(pappbp)=f(pap)f(pbp)=f(a)f(b)$ because $p\Nc$ is
1-dimensional.  

If $f$ is multiplicative and normal, then $(1-p)\Nc$ is  
an ideal and the kernel of $f$, so $1-p$ is central.  This finishes part a.

By part a, if $f$ is not singular (i.e. $f$ is normal), then no approximate unit
of
$K_f$ could converge weak* to 1 because each element of $K_f$ vanishes on $f$,
and hence
$f(1)=0$, contradicting the assumption that $f$ is a state.

Now suppose that $f$ is singular and some increasing, positive approximate unit
(of $K_f$)
$a_\alpha \uparrow r \ne 1$ in the weak* topology of $\Nc$.  By \cite{tak2},
there is a projection
$p
\le 1-r$ such that $f(p)=0$.  This contradicts the fact that $\{a_\alpha\}$ is an
approximate unit for $K_f$ since $a_\alpha p=0$ for all $\alpha$.

\end{proof}

\begin{lemma}\label{LM}
A  MASA $\Ac$ of a C*-algebra $\Bc$ contains no
minimal projections that are not minimal in all of $\Bc$.
\end{lemma}
\begin{proof}
If $p$ is a minimal projection of $\Ac$ and $p$ is not minimal in $\Bc$, then
$pBp$ contains a projection $q$ that is not $p$ or $0$.  However, for any $a
\in A$, $qa-aq=qpa-apq=qpap=papq=q(\lambda p)-(\lambda p)q=0$ because
minimality of $p$ in $\Ac$ implies that $pap$ is a multiple of $p$ for every $a
\in A$.

\end{proof}

\begin{lemma}\label{LM}
A pure state $f$ of a von 
Neumann algebra $\Mx$ is {\it classically normal} if $f$ is normal on any
abelian von Neumann subalgebra of $\Mx$  on which $f$ is
multiplicative. 

\end{lemma}
\begin{proof}
This follows immediately from \cite{tak} Cor. III.3.11 and the remark following. 
\end{proof}

\noindent {\bf NOTATION.} Now we fix a factor $\Nc$ of type $I_\infty, II$ or
$III$ on separable Hilbert space and WLOG assume that a type $I_\infty$ factor is
all of $\Bc(H)$.  We let $\Cc(\Nc)$ denote
the set of all von Neumann subalgebras of $\Nc$.

\begin{lemma}\label{LM} Suppose
that $\{p_n\}$ is a decreasing sequence of projections in $\Nc$ that
converges to 0 in the weak-operator topology. If $\Cc \in C(\Nc)$, then 
there is a singular pure state
$f$ on $\Nc$ and a projection $q$ in
$\Cc$ such that $f(p_n)=1 \ \forall n$ and either

1.  $f(q) \in (0,1)$, or

2.  $q$ is minimal in $\Cc$, $q$ is a central projection in $\Cc$, and $f(q)=1$.

\end{lemma}

\begin{proof}
Working in $N^{**}$ (see \cite{gert}, sect.
3.8), let $p = \lim\ p_n$. By Corollary 4.5.13 of \cite{gert} 
there exists at least
one pure state $f$ of $\Nc$ such that $f(p) =1$.  Any such $f$ is singular
since $p_n \downarrow 0$ in the weak* topology of $\Nc$.  We need to  find such
an
$f$ and a projection $q \in \Cc$ such that 1 or 2 of the Lemma holds.  

Suppose that 1 does not hold for any pure state $f$ of
$\Nc$ with
$f(p)=1$ and projection $q \in \Cc$.  By Lemma 0.1, any pure state $f$ of
$\Nc$ with
$f(p)=1$ is multiplicative on $\Cc$. (Using \cite{kr}, 
Theorem 6.5.2, and a bit more
argument, we can assume that
$\Cc$ is abelian, since only the abelian direct summand of a von Neumann algebra
can support a nonzero multiplicative linear functional.  This is not required for
the proof, but it does serve to clarify matters.)  We now proceed by cases to
reach a contradiction or conclude that condition 2 of the Lemma holds.

{\bf Case 1:}  There is a state $g$ of $\Cc$ such
that, if $h$ is any pure state of $\Nc$ such that $h(p)=1$, then $h|_{\Cc} =
g$.

Let $Q$ denote the set of all states $e$ of $\Nc$ such that  $e(p)=1$.  Then
$Q$ is convex and weak* compact face of the state space of $\Nc$ by
\cite{face}, Theorem 2.10.  By the Krein-Milman Theorem,
\cite{rud},Theorem 3.21, $Q$ is the weak* closed convex hull of its extreme
points, and the extreme points of
$Q$ are pure states of $\Nc$.  Since any of the extreme points of $Q$ restrict
to $g$ on $\Cc$, the same must be true of all the states of $Q$. As noted above,
$g$ is multiplicative on $\Cc$.

There are two subcases:

{\bf Subcase a.}  
$g$ is normal on $\Cc$. 
Choose any pure state $f$ of $\Nc$ such that $f(p)=1$. Then $q$ exists
satisfying condition 2 of the lemma by Lemma 0.2a.

Let $\Cc_g = \{a \in \Cc : g(a) = 0\}$, and let
$\{r_\alpha\}$ be an approximate unit in
$\Cc_g$ with

{\bf Subcase b.} $g$ is singular on $\Cc$. Let $f$ be a pure state
 of $\Nc$ such
that $f(p)=1$. Since $f|{\Cc}=g$ is singular, $f$ must also be singular.  By
part b of Lemma 0.2, if $\{r_\alpha\}$ is an increasing approximate unit for 
$K_f=\{a\in\Nc : f(a^*a+aa^*) = 0\}$ with $r_\alpha
\uparrow r \in \Nc^{**}$, then $r_\alpha \to 1$ in the weak*
topology of $\Nc$.  Then $r$ is an open projection,
hence regular in $\Nc$ by \cite{sw}, Prop.
II.14 since $\Nc$ is a von Neumann
algebra.
 Also, if $\overline r$ denotes the closure of $r$ in $\Nc^{**}$, then
${\overline r} \in \Nc$ by \cite{left}, Theorem II.1, so $\overline
r=1$. 

 Thus
$\|p_nrp_n\|=\|rp_n\|^2=1$ for all $n$ by regularity of $r$.  Consequently
 by Corollary 4.5.13 of \cite {gert}, there
exist pure states  $\{f_n\}$ on $\Nc$ such that
$|f_n(p_nr p_n)|>1-1/n$ and $f_n(p_n)=1$ for all $n$.  Set $g_n = \ f_n|_\Cc$. 
Since any limit point $h$ of $\{f_n\}$ in
$N^*$ must satisfy $h(p_n) = 1$ for all $n$,
then
$h(p) = 1$.  Consequently, by the paragraph following the
Case 1 assumption,
$h|_C = g$ and  
$f_n|_C \to g$ in $\Cc^*$ (weak* topology).  However,
since
$f_n(r)> .5$ for all large $n$ and $g(r) = 0$,  we
contradict \cite{seq}, Theorem 4. Thus subcase b leads to a contradiction, so
subcase a must hold for Case 1.

{\bf Case 2:}  The remaining possibility is that there are two pure
states $f,g$ of $\Nc$ such that $f|_C \ne g|_\Cc$ and $f(p)=g(p)=1$.   As noted
in the paragraph above the Case 1 statement, $f$ and $g$ are multiplicative on
$\Cc$, hence there exists a projection $q$ in $\Cc$ such
that $f(q) = 1, g(q) = 0$. 
Since $p$ commutes with $\Cc$ by Corollary 4.5.13 of \cite{gert} and the
assumption that property 1 of the Lemma is false,
$q-pq = 1-(\sup(p,(1-q))$, is open by Corolary II.7 of \cite{sw}. Since $q-qp_nq
\uparrow q-qp$, the operator $b=\sum_1^\infty 2^{-n}(q-qp_nq)$ is strictly
positive in the hereditary C*-subalgebra  
$$D=\{d \in N:\|(q-qp_nq)d\|+ \|d(q-qp_nq)\|\to 0\}.$$  Thus the spectral
projections
$r_n=q-q\chi_{(1/n,\|b\|]}(b) \downarrow q-(q-qp) = qp$.
Similarly we get a sequence
$\{s_n\}$ of projections such that  $1-q \ge s_n \downarrow (1-q)p$.
 Since
$p_n\to 0$ in the weak operator topology, we can pass to a subsequence and
assume that
$r_n-r_{n+1}$ and  $s_n-s_{n+1}$ are all nonzero.  To reach a contradiction of
the assumption that Case 1 fails, it suffices to find an irreducible
representation $\pi$ of $\Nc$ such that $\pi^{**}(pq)\ne 0 \ne \pi^{**}(p(1-q))$.

N.B.  Up to this point in the proof, $\Nc$ could be any von Neumann
algebra.  We now break the proof into subcases to handle the different types
of factors.

{\bf Type $III$ subcase:}
Set $e_n = (r_n - r_{n+1}), e'_n= (s_n-s_{n+1})$.  
Set $e=\sum_{n=1}^\infty e_n, e'=\sum_{n=1}^\infty
e'_n$. Choose partial isometries $v_n$ such that $v_n^*v_n = e_n',
v_nv_n^*=e_n$, and set $v=\sum_1^\infty v_n$ (any two non-trivial
projections are equivalent because $\Nc$ is a type $III$ factor and separably
represented).  By  Corollary 4.5.13
of
\cite{gert}, we can choose an irreducible reepresentation $\pi$ of $\Nc$ such that
$$\pi^{**}(pq)=\lim\pi(r_n)\ge\lim_{k\to\infty}\pi(\sum_{n=k}^\infty e_n) \ne 0
$$ because $\lim_{k\to\infty}\sum_{n=k}^\infty e_n \ne 0$.

Hence
$$\|\pi^{**}(p(1-q))\|=\|\lim_{k \to \infty}\pi(s_k)\|\ge\|\lim_{k\to\infty}
\pi(\sum_{n=k}^\infty e'_n)\|\ge\|\lim_{k\to\infty}
\pi(v^*)\pi(\sum_{n=k}^\infty e_n)\pi(v)\|$$
$$=\|\pi(v^*)(\lim_{k\to\infty}\pi(\sum_{n=k}^\infty e_n))\pi(v)\|$$
$$\ge \|\pi(v)\pi(v^*)(\lim_{k\to\infty}\pi(\sum_{n=k}^\infty
e_n))\pi(v)\pi(v^*)\| =\|\pi(e)(\lim_{k\to\infty}\pi(\sum_{n=k}^\infty
e_n))\pi(e)\|$$
$$=\|(\lim_{k\to\infty}\pi(\sum_{n=k}^\infty e_n))\|\ne 0$$  as
was to be shown.

{\bf Type $II$ subcase:}  Let $\tau$ denote a normal, semi-finite trace on
$\Nc$.
  We use the same idea as in the type $III$ case 
except we use the continuity of the trace to choose 
$e_n \le (r_n -
r_{n+1}), e'_n\le (s_n-s_{n+1})$ such that 
$0 < \tau(e_n) = \tau(e'_n) \le 2^{-n}$ for all $n$. 
Choose partial isometries $v_n$ and define $v$ as in the type $III$ case
(because any two projections with the same finite trace are equivalent). 
The same contradiction then arises
as in the type
$III$ subcase.

{\bf Type $I_\infty$ subcase:} Again we mimic the type $III$ case 
with the following exception. Since $r_n \ne 0$ and $s_n \ne 0$,.  We can then
choose rank 1 projections
$e_n\le (r_n - r_{n+1}), e'_n\le (s_n-s_{n+1})$.  
Choose partial isometries $v_n$ and define $v$ as in the type $III$ case.
The contradiction
then follows as in the Type $III$ subcase.

\end{proof}

\begin{lemma}\label{LM} Let $\Ac$ denote a separable C*-subalgebra of
$\Nc$ with the same unit, and assume that $\Ac$ contains the compact operators
in the type $I_\infty$ case. Let $\Cc \in C(\Nc)$.  Suppose that
$h$ is a pure state of
$\Ac$ that annihilates the conpact operators in $\Ac$. Then
there is a singular pure state $f$ on
$\Nc$ that extends $h$ and a projection $q \in \Cc$ such that either 

1.  $f(q) \in (0,1)$, or

2.  $q$ is minimal in $\Cc$, $q$ is a central projection in $\Cc$, and $f(q)=1$.

\end{lemma}

\begin{proof}
Let $\Ac_h = \{a \in A : h(a^*a+aa^*) = 0\}$.  Since $\Ac_h$ is separable, then
it contains a completely positive  element $a$ of norm 1. 

There are two cases: 

Case 1.  $0$ is an isolated point in the spectrum of $a$. 
Applying the functional calculus to $a$, we can get a projection $r \in K_h$
that acts as a unit for $K_h$.  $1-r$ must be of infinite rank in $\Nc$ because
by assumption $K_h$ contains all the finite rank projections in $\Nc$.  Let $g$
be any singular pure state of $\Nc$ such that $g(1-r)=1$.  Since $g$ is singular
and
$\Nc$ acts on a separable Hilbert space, by
\cite{tak2} there is a decreasing sequence $\{p_n\}$ of projections in $\Nc$ such
that $p_1 \le 1-r$ and $p_n\downarrow 0$ in the weak* toplogy of $\Nc$ and
$g(p_n)=1$ for all $n$. Let $B$ be the C*-algebra generated by $A$ and
$\{p_n\}$. Clearly $h = (1-r)h(1-r)$ and $g = (1-r)g(1-r)$, and on
$(1-r)B(1-r)\cap \Ac = \{\lambda (1-r):\lambda \in $C$\}, h=g$.  Thus $g$ extends
$h$.  Further, any pure state $f$ of $\Nc$ such that $f(p_n)=1$ for all $n$ will
extend $g$ (and hence extend $h$ also).  Lemma 0.5 now gives the desired $f$ and
$q \in \Cc$. 

Case 2.  If 0 is not an isolated point in the spectrum of $a$, then set
$p_n=\chi_{[0,1/n]}(a)$.  
Let $B$ be
the C*-algebra generated by $\Ac$ and $\{p_n\}$.  Let $g$ be any extension of $h$
to a pure state of $B$.  Then $g(1-p_n) \le g(na)=h(na)=nh(a)=0$, so $g(p_n)=1$
for all $n$.   We need only show that any pure state $f$ of $\Nc$ such that
$f(p_n)=1$  for all $n$ will
extend $g$; then Lemma 0.5 gives the desired $f$ and
$q \in \Cc$.

Let $f$ be any pure state of $\Nc$ such that $f(p_n)=1$  for all $n$.  Let $c
\in \Ac$.  Let $p = \lim p_n$ in $\Bc^{**}\subset \Nc^{**}$.  Then 
$\{1-p_n\}$ is an approximate unit for $\{d \in \Bc : g(dd^*+d^*d)=0$, so by
\cite{excising} Proposition 2.2, p is a minimal projection in $\Bc^{**}$ and
$g=pgp$.  But $f(p) = 1$ also, so $f|_{\Bc}=g|_{\Bc}$, and the lemma follows.

\end{proof}

\begin{theo}\label{TH}
Assume the continuum hypothesis. There is a classically normal, singular pure
state $f$ on $\Nc$.  In particular, $f$ is not multiplicative on any MASA of
$\Nc$
\end{theo}

\begin{proof}
 Let
$(a_\alpha)$,
$\alpha <
\aleph_1$, enumerate the elements of
$\Nc$. Since every von Neumann subalgebra of $\Nc$ is countably generated, a
simple cardinality argument shows that the cardinality of $\Cc(\Nc)$ is
$\aleph_1$. Let
$(\Mx_\alpha)$,
$\alpha <
\aleph_1$, enumerate $\Cc(\Nc)$.

We now inductively construct a nested transfinite sequence of unital
separable C*-subalgebras $\Ac_\alpha$ of $\Nc$ together with pure
states $f_\alpha$ on $\Ac_\alpha$ such that for all $\alpha < \aleph_1$
\begin{enumerate}
\item
$a_\alpha \in \Ac_{\alpha+1}$
\item
if $\beta < \alpha$ then $f_\alpha$ restricted to $\Ac_\beta$ equals
$f_\beta$
\item
$\Ac_{\alpha+1}$ contains a projection $q_\alpha \in \Mx_\alpha$ such that either
$0 < f_{\alpha+1}(q_\alpha) < 1$ or else $q_\alpha$ is minimal
and central in $\Ac_\alpha$
with
$f_{\alpha+1}(q_\alpha)=1$.
\end{enumerate}
Begin by letting $\Ac_0$ be any separable C*-subalgebra of $\Nc$
that is unital (and contains $\Kc(H)$ when $\Nc$ is type $I_\infty$) and let
$f_0$ be any pure state on
$\Ac_0$ (that annihilates $\Kc(H)$ when $\Nc$ is type $I_\infty$). At
successor stages, use the last lemma to find a projection $q_\alpha \in
\Mx_\alpha$ and a pure state
$g$ on $\Nc$ such that $g|_{\Ac_\alpha} = f_\alpha$ and either
$0 < g(q_\alpha) < 1$ or else $q_\alpha$ is minimal and central in $\Ac_\alpha$
and
$g(q_{\alpha})=1$. By
\cite{cn}, Lemma 4, there is a separable C*-algebra $\Ac_{\alpha+1} \subseteq
\Nc$ which contains
$\Ac_\alpha$,
$a_\alpha$, and $q_\alpha$ and such that the restriction $f_{\alpha+1}$
of $g$ to $\Ac_{\alpha+1}$ is pure. To see this, write $\Nc$ as the
union of a continuous nested transfinite sequence of separable C*-algebras
$\Bc_\gamma$ such that $\Bc_0$ is the C*-algebra generated by $\Ac_\alpha$,
$a_\alpha$, and $q_\alpha$. The cited lemma guarantees that the restriction
of $g$ to some $\Bc_\gamma$ will be pure. Thus the construction may proceed.
At limit ordinals $\alpha$, let $\Ac_\alpha$ be the closure of
$\bigcup_{\beta<\alpha} \Ac_\beta$. The state $f_\alpha$ is determined
by the condition $f_\alpha|_{\Ac_\beta} = f_\beta$, and it is easy to
see that $f_\alpha$ must be pure. (If $g_1$ and $g_2$ are states on
$\Ac_\alpha$ such that $f_\alpha = (g_1 + g_2)/2$, then for all
$\beta < \alpha$ purity of $f_\beta$ implies that $g_1$ and $g_2$
agree when restricted to $\Ac_\beta$; thus $g_1 = g_2$.) This
completes the description of the construction.

Now define a state $f$ on $\Nc$ by letting $f|_{\Ac_\alpha} = f_\alpha$.
By the reasoning used immediately above, $f$ is pure, so $f$ is a classically
normal, singular pure state.

If $\Ac$ is a MASA of $\Nc$, then in the type II or II cases, there are no
minimal projections in $\Nc$, hence none in $\Ac$, so $f$ can't be
multiplicative on $\Ac$.  In the type I case, the only minimal projections are
the rank one projections.  Since $f$ is singular, it must vanish on rank one
projections, hence it can't be multiplicative (and hence normal because it is
classically normal) on $\Ac$.
\end{proof}

\begin{coro}\label{TH}
Assuming the continuum hypothesis, a separably represented von
Neumann algebra $M$
has classically normal, singular pure states iff there is a central projection
$p \in M$ such that $pMp$ is a factor of type $I_\infty, II$, or $III$.
\end{coro}

\begin{proof}
The implication $\leftarrow$ was essentially proved in the last Theorem since any
pure state on $pMp$ will have unique pure extension to $M$.  

For the other direction, suppose that no such
projection $p$ exists, so that the center of $M$ is infinite dimensional. 
Suppose that $f$ is a classically normal singular pure state of $\Nc$. Since any
pure state of a C*-algebra is always multiplicative on the center of the
algebra, $f$ must be normal on the center of $\Nc$ by the definition of
classically normal.  Thus by Lemma 0.2a there is a minimal projection $p$ in the
center of
$\Nc$ such that $f(p)=1$.  Since $p$ is minimal in the center, $p\Nc p$ is a
factor.  By the assumption of the Corollary, $p\Nc p$ must be a factor of
type $I_n$ for $n < \infty$, i.e. $p$ must be a finite rank projection in $\Nc$.
However, a singular state of
$\Nc$ must vanish on such projections by Takesaki's singularity cirterion
\cite{tak2}.  Since $p$ is a finite sum of minimal projections in $\Nc$ on
which $f$ must vanish, $f(p)=0$. This contradiction completes the proof.

\end{proof}

REMARK: Since there is a choice to make at each of the $\aleph_1$ steps of the
proof, assuming CH, the methods of the last theorem will produce
$2^{\aleph_1}$ classically normal pure states on $\Nc$.  Since $\Nc$ has
cardinality
$\aleph_1$, the totality of states of $\Nc$ must be of cardinality
$2^{\aleph_1}$.  Since any MASA of an infinite factor has (under CH)
$2^{\aleph_1}$ distinct singular pure states, each of which has an extension to
a pure state of
$\Nc$, there must be $2^{\aleph_1}$ pure states that are not classically normal. 
\bigskip

\begin{coro}\label{TH}
Let $f$ be a classically normal, singular pure state of $\Nc$.  Let
$\Ac = \{a \in N: f(a^*a+aa^*)=0\}$.  Then A does not have an abelian approximate
unit.
\end{coro}

\begin{proof}
Suppose the Corollary is false. I.e.  $f$ is a classically normal, singular pure
state of $\Nc$ such that $\Ac = \{a \in N: f(a^*a+aa^*)=0$ does  have
an abelian approximate unit $\{a_\alpha\}$.  Let $\Bc$ be a MASA of $\Nc$ that
contains $\{a_\alpha\}$.  Then $\Bc$ is
an abelian von Neumann subalgebra of $\Nc$ that contains a decreasing excising
net
$\{1-a_\alpha\}$ for
$f|_B$ by
\cite{excising}, Prop. 2.3, and $f|_{\Bc}$ must be pure, hence multiplicative.
Since
$f$ is classically normal, there is a minimal projection $q \in \Bc$ such that
$f(q)=1$.  By Lemma 0.3, $q$ is minimal in $\Nc$ also, so $f$ is not
singular by Lemma 0.2, a contradiction. 
\end{proof}

Not all questions of this general type are resolvable by the methods above.  For
instance, if for each natural number $n$, $H_n$ is a Hilbert space of dimension
$n$, and if $M = \sum_1^\infty \oplus B(H_n)$, then $M$ does not have any
singular, classically normal pure states.  However, our methods don't say whether
or not $M$ has a pure state that is not multiplicative on any MASA.

\bigskip
We conclude the paper by mentioning an example that shows that there is more to
the existence of classically normal pure states than substantial
non-commutativity and nonseparability of the underlying algebra.
\bigskip

NOTATION:  $F_R$ is the free group on card($R)=\aleph_1$
 generators and
$\Cc^*(F_R)$ is the corresponding reduced group C*-algebra.

\bigskip
This example is discussed in \cite
{pop}, Cor. 6.7, where it is shown that $\Cc^*(F_R)$ is nonseparable, but that
every abelian subalgebra is separable.  However, unlike the situation
described in Corollary 0.9, we show in \cite{cn3} that, if $f$ is a pure state on
$\Cc^*(F_R)$, then 
$\Ac = \{a \in C^*(F_R): f(a^*a+aa^*)=0\}$ contains a sequential abelian
approximate unit.

\end{document}